\documentstyle[12pt]{article}
\begin{document} 

\Large 
{\bf 
Riccati equations and convolution formulas for functions of
Rayleigh type}\\

\large  
\begin{quote}
Dharma P. Gupta and
Martin E. Muldoon\footnote{e-mail address: {\tt muldoon@yorku.ca}.
Research supported by grants from the
Natural Sciences and Engineering Research Council, Canada}
\normalsize
\\ Department of Mathematics \& Statistics,
York University\\
Toronto ON M3J 1P3, Canada\\
\bigskip

Revised \today 
\bigskip

{\bf Abstract.} N. Kishore, {\em Proc. Amer. Math. Soc.} {\bf 14}
527 (1963), considered the Rayleigh functions
$\sigma_n(\nu) = \sum_{k=1}^\infty j_{\nu k}^{-2n},\;n=1,2,\dots$,
where $\pm j_{\nu k}$ are the (non-zero) zeros of the Bessel
function $J_\nu(z)$ and provided a convolution type sum formula
for finding
$\sigma_n$ in terms of $\sigma_1,\dots,\sigma_{n-1}$.  His main
tool was the recurrence relation for Bessel functions. Here we
extend this result to a larger class of functions by using Riccati
differential equations.  We get new results 
for the zeros of certain combinations of Bessel functions and their
first and second derivatives as well as recovering some results of
Buchholz for zeros of confluent hypergeometric functions.
\end{quote}

\normalsize
\section{Introduction}
The 
Rayleigh functions are defined, e.g., in  \cite[p. 502]{wat}, by
the formula
\begin{equation}
\sigma_n(\nu) = \sum_{k=1}^\infty j_{\nu k}^{-2n}, \;n=1,2,\dots,
\end{equation}
where  $\pm j_{\nu k}$ are the zeros of the Bessel function
\begin{equation}
J_{\nu}(z) = \sum_{n=0}^{\infty}\frac{(-1)^n
(z/2)^{2n+\nu}}{n!\Gamma(\nu+n+1)}.
\label{series}
\end{equation}
They form the basis of an old method due to Euler, Rayleigh and
others for evaluating  the zeros.
For example, in case $\nu > -1$, the inequalities 
$$
[\sigma_n(\nu)]^{-1/n} < j_{\nu 1}^2 < \sigma_n(\nu)/
\sigma_{n+1}(\nu),\;\; \;n=1,2,\dots
$$
provide infinite sequences
of successively improving upper and lower bounds for $j_{\nu 1}^2$.
Several authors have considered the question of finding ``sum
rules" or formulas for
$\sigma_n(\nu)$. By a method originating with Euler (see \cite[p.
500, {\em ff.}]{wat} for details; various ramifications 
were considered recently in \cite{im}), we can find
all the $\sigma_n(\nu)$ in terms of the coefficients in the series 
(\ref{series}).
If we want to deal (as in \cite{obi}) with properties of the 
$\sigma_n(\nu)$ as functions of $\nu$, there is a useful
compact convolution formula 
due to Kishore \cite{kish} \begin{equation}
 \sigma_n(\nu) = \frac{1}{\nu +n}\sum_{k=1}^{n-1}\sigma_k(\nu)
\sigma_{n-k}(\nu),\label{kish1}
\end{equation}
from which the $\sigma_n(\nu)$ may be found successively, starting
from
\begin{equation}
\sigma_1 =1/[4(\nu +1)].
\label{sigma1}
\end{equation}

The question arises whether there are Kishore-type
formulas for sums of
zeros of other special functions such as the
first and second derivatives of the Bessel function.  In
\cite{raza} there is 
a variant of this result for the zeros of the more general function
\begin{equation}
N_\nu(z) = az^2J_\nu^{''}(z) + bzJ_\nu^{'}(z) + cJ_\nu (z) 
\label{defn}
\end{equation}
considered by Mercer \cite{mercer}.  
The result of \cite{raza} gave a method of finding the reciprocal
power sums
\begin{equation}
\tau_n(\nu) = \sum_{k=1}^\infty x_{\nu k}^{-2n}, \;n=1,2,\dots.
\end{equation}
where $x_{\nu k}$ are the zeros of the function $N_\nu(z)$.
The main result of \cite{raza} expressed $\tau_n$ in terms of 
$\tau_k,\;k=1,\dots,n-1$ {\em and} $\sigma_k,\;k=1,\dots,n$.
It seems desirable to express 
$\tau_n$ in terms of 
$\tau_k,\;k=1,\dots,n-1$ only.  We do this here by using the
Riccati equation satisfied by $z^{-\nu/2}N_\nu(x^{1/2})$.  We
record also the
second order linear differential equations satisfied by $N_\nu(z)$
and by $z^{-\nu/2}N_\nu(z^{1/2})$ since these do not seem to appear
in the literature and may prove useful for other purposes. 

In \S4, we apply the same method to get power sums for zeros of
confluent hypergeometric functions.

\section{Differential equations for functions related to Bessel
functions}

The Bessel function $ y = J_\nu(z)$ 
satisfies the differential equation
\begin{equation}
z^2y'' +zy' + (z^2 - \nu^2)y = 0. \label{bde}
\end{equation}
and the function $y = zJ_\nu^{'}(z) + cJ_\nu (z) $
satisfies \cite[p. 13]{erdelyi} the differential equation
\begin{eqnarray*}
z^2(z^2 - \nu^2+c^2)y'' &-& z(z^2 + \nu^2-c^2)y' \\
&+& [(z^2 - \nu^2)^2 + 2cz^2 +c^2(z^2 - \nu^2)]y = 0. 
\end{eqnarray*}
Here we record the more general second order linear differential
equation satisfied by the function
\begin{equation}
y = N_\nu(z) = az^2J_\nu^{''}(z) + bzJ_\nu^{'}(z) + cJ_\nu (z) .
\label{defn1}
\end{equation}
It is 
\begin{equation}
z^2y''+ A(z)zy' + [B(z) + z^2 - \nu^2]y = 0, \label{neweq}
\end{equation}
where
$$ A(z) = \frac
{-3a^2z^4 + pz^2 + q}
{a^2z^4 - pz^2 + q},
$$
$$ B(z) = \frac
{2a(a + b)z^4 + 2rz^2}
{a^2z^4 - pz^2 + q},
$$
 with
$$ p = 2a(a\nu^2+c) +(a^2 - b^2),$$
$$ q = (a\nu^2+c)^2 -\nu^2(a - b)^2,$$
and $$r = a\nu^2(3a-b) + c(a+b). $$
We found the equation (\ref{neweq}) by repeated use of 
\begin{equation}
zJ_\nu'(z) = \nu J_\nu(z) -zJ_{\nu+1}(z)
\end{equation}
to express the derivatives $J_\nu^{(n)}(z),\;
n=1,2\dots$ in terms of
$J_\nu(z),\;J_{\nu+1}(z)$ and discovered an appropriate vanishing
linear combination of $N_\nu(z),\;N_\nu'(z),\;N_\nu''(z)$.  Of
course, once (\ref{neweq}) is known, it is easy to verify that
$N_\nu(z)$, given by  (\ref{defn1}), satisfies it.

It is convenient to consider the function
\begin{equation}
y_\nu(z) = z^{-\nu/2}N_\nu(z^{1/2}) 
\end{equation}
where we choose that branch of $z^{1/2}$ which is positive for $z
>0$. Using (\ref{neweq}), we find that the function $ y_\nu(z) $
satisfies
\begin{equation}
4t^2\frac{d^2y}{dt^2} + [4\nu + 2 + 2A(t^{1/2})]t\frac{dy}{dt} 
+[t -\nu +\nu A(t^{1/2})+ B(t^{1/2})]y =0, \label{neweqq}
\end{equation}
It is well known that if $y$ satisfies
\begin{equation}
y'' + P(t)y' + Q(t)y = 0,
\end{equation}
then
$u = y'/y$ satisfies
the Riccati equation
\begin{equation}
\frac{du}{dt} +P(t)u + Q(t) + u^2 = 0. \label{riccmer}
\end{equation}
Applying this to (\ref{neweqq}), we find that, with $y_\nu(z)$
given by
(\ref{defn1}), $u = y_\nu'(z)/y_\nu(z)$ satisfies
\begin{eqnarray}
4t(a^2t^2  - pt+q)\left[\frac{du}{dt} \right.&+& \left. u^2 \right]
+
4[a^2(\nu-1)t^2 -\nu pt +q(\nu +1)]u +
\nonumber\\
&+& a^2t^2 +[p+4a^2\nu -2a(a+b)]t +2\nu p +q + 2r
= 0. \label{riccb}
\end{eqnarray}
\section{Functions of Rayleigh type}
The even entire function $z^{-\nu}N_\nu(z)$ has an infinite set of
zeros
$\pm t_n,\;n=1,2,\dots$
with
$$\sum |t_k^{-2}| < \infty,$$
so the zeros of $y_\nu(z)$ are $\zeta_k = t_k^2,$
with
$$\sum |\zeta_k^{-1}| < \infty.$$
Thus
\begin{equation}
y_\nu(z) = z^{-\nu/2}N_\nu(z^{1/2}) = \frac{a\nu^2 + c
+(b-a)\nu}{2^\nu
\Gamma(\nu +1)}
\prod_{k=1}^\infty \left( 1 - \frac{z}{\zeta_k} \right).
\label{prod}
\end{equation}
The constant multiplicative factor is got from the series
(\ref{series}). The validity of this infinite product expansion
follows from facts on entire functions of finite order \cite[Ch.
8]{titchmarsh}. 

We may differentiate (\ref{prod}) logarithmically \cite{knopp}, to
get
$$
\frac{y_\nu'(z)}{y_\nu(z)} = - \sum_{k=1}^\infty
\frac{1/\zeta_k}{1- z/\zeta_k} =
- \sum_{k=1}^\infty \frac{1}{\zeta_k}
\sum_{n=0}^\infty \frac{z^n}{\zeta_k^n}.
$$
This gives
$$
2z\frac{y_\nu'(z)}{y_\nu(z)} =
 -
2 \sum_{k=1}^\infty \sum_{n=1}^\infty z^n/\zeta_k^n.
$$
But we may interchange the orders of summation here (since the
iterated series converges absolutely) to get
\begin{equation}
2z\frac{y_\nu'(z)}{y_\nu(z)} = -
2 \sum_{n=1}^\infty z^n\sum_{k=1}^\infty \zeta_k^{-n}
= - 2\sum_{n=1}^\infty
\tau_nz^n, \label{19}
\end{equation}
where
\begin{equation}
\tau_n = \sum_{k=1}^\infty \zeta_k^{-n}.
\end{equation}
Using 
$$
u = -\sum_{k=0}^\infty \tau_{k+1}z^k, $$
we get
$$
u^2 = \sum_{k=2}^\infty \left[ \sum_{m=1}^{k-
1}\tau_m\tau_{k-m} \right]z^{k-2}. $$
Substituting in (\ref{riccb}), and comparing coefficients of powers
of
$z$, we get
\begin{equation}
\tau_1 = \frac{2\nu p+q+2r}{4q(\nu+1)},$$
$$ 4q(\nu+2)\tau_2 = 4q\tau_1^2+4 \nu p \tau_1 - p - 4a^2\nu
+ 2a(a+b), \label{res1}
\end{equation}
\begin{equation} 4q(\nu+3)\tau_3 = 4p(\nu+1)\tau_2 - 4a^2(\nu-1)
\tau_1 +a^2
+8q\tau_1\tau_2 - 4p\tau_1^2, \label{res2} \end{equation}
and, for $k \ge 3$,
\begin{eqnarray}
q(k + \nu +1) \tau_{k+1}&=&p(k + \nu -1) \tau_{k} -a^2(k + \nu -3)
\tau_{k-1} \nonumber\\
& &
+q \sum_{m=1}^{k} \tau_m \tau_{k-m+1}
-p\sum_{m=1}^{k-1} \tau_m \tau_{k-m}
+a^2 \sum_{m=1}^{k-2} \tau_m \tau_{k-m-1} \label{res3}
\end{eqnarray}

In the special case $a=b=0,\;c=1$ (and hence $p=0,\;q=1,\;r = 0$),
where we are dealing with the
zeros of the Bessel function, these reduce, as they should, to
(\ref{sigma1}) and the convolution formula (\ref{kish1}) for
$\sigma_n,\;n=2,3,\dots$.

In the special case
$a=c=0, b=1$ (and hence $ p = -1, q = -\nu^2, r = 0$), we are
dealing with the
non-trivial
zeros of the function $J_\nu'(z)$; (\ref{res1}), (\ref{res2}) and
(\ref{res3})
become
\begin{equation}
\tau_1 = \frac
{ \nu +2 }
{4(\nu+1)\nu}   \label{res12}
\end{equation}
$$ 
\tau_2 = \frac{- 4\nu^2\tau_1^2 -4 \nu  \tau_1 +1}{-
4\nu^2(\nu+2)},$$
$$ \nu^2(\nu+3)\tau_3 = (\nu+1)\tau_2  +2\nu^2\tau_1\tau_2 -
\tau_1^2,$$
and for $k \ge 3$,
\begin{eqnarray}
-\nu^2(k + \nu +1) \tau_{k+1}&=&-(k + \nu -1) \tau_{k}  \nonumber\\
& &
-\nu^2 \sum_{m=1}^{k-1} \tau_m \tau_{k-m+1}
+\sum_{m=1}^{k-2} \tau_m \tau_{k-m}
\end{eqnarray}
In particular, these lead to
\begin{equation}
\tau_2 = \sum_{k=1}^\infty [j_{\nu k}']^{-4} = 
\frac{1}{16} \frac{\nu^2+8\nu+8}{\nu^2(\nu+1)^2(\nu+2)},
\end{equation}
\begin{equation}
\tau_3 = \sum_{k=1}^\infty [j_{\nu k}']^{-6} = 
\frac{1}{32}\frac{\nu^3+16\nu^2+38\nu+24}{\nu^3(\nu+1)^3(\nu+2)
(\nu+3)},
\end{equation}
the same results as are obtained by the power series method
in \cite{im}.
\section{Confluent Hypergeometric Functions}
Buchholz \cite{buchholz} studied the nontrivial zeros $a_\lambda$
of
the function
\begin{equation}
M_{\kappa, \mu/2}(z) = \frac{z^{b/2}e^{-z/2}}{\Gamma(1 + \mu)} \;
 _1F_1(a;b;z)
\end{equation}
and showed that these zeros are all simple and that there are
infinitely many of them in the case where $a \ne -n$.
He considered 
$$ S_p =
\sum_{\lambda = 1}^\infty a_\lambda^{-p}, $$ 
and showed that it converges
for all $p > 1$ but that it is divergent for
$p \le 1$. 

He also gave explicit formulas for $S_2,\; \dots\;S_{6}$ and a
method
(far from explicit) for expressing $S_{k+1}$ as a linear
combination of $S_2,\; \dots\;S_{k-1}$. In (\ref{kishchf}) below we
give a
convolution formula for this task.

The function $w = \;_1F_1(a;b;z)$ satisfies
\begin{equation}
zw'' + (b-z)w' -aw = 0
\end{equation}
so $u = w'/w$ satisfies the Riccati equation
\begin{equation}
zu'+ (b-z)u -a +zu^2 = 0. \label{riccc}
\end{equation}
From the Weierstrass product representation theorem, we get
\begin{equation}
w = e^{az/b} \prod_{k=1}^\infty \left(1 - \frac{z}{z_k}
\right)e^{z/z_k}. \label{prodh}
\end{equation}
Differentiating (\ref{prodh}) logarithmically \cite{knopp}, 
\begin{eqnarray}
u(z) =\frac{w'(z)}{w(z)} &= & \frac{a}{b} - \sum_{k=1}^\infty
\left[
\frac{1/z_k}{1- z/z_k} - \frac{1}{z_k} \right] \nonumber \\&=&
\frac{a}{b} - \sum_{k=1}^\infty \frac{1}{z_k}
\left\{ \left[ 1 - \frac{z}{z_k} \right]^{-1} -1 \right\}\nonumber
\\
&=&
\frac{a}{b} - \sum_{k=1}^\infty S_{k+1}z^k, \label{u}
\end{eqnarray}
where the interchange of orders of summation here is justified
by the absolute convergence of the 
iterated series. 
From this we have
\begin{equation}
zu'(z) = 
- \sum_{k=1}^\infty kS_{k+1}z^k,
\end{equation}
and
\begin{equation}
[u(z)]^2 = (a/b)^2 - 2(a/b)
\sum_{k=1}^\infty S_{k+1}z^k
+ \sum_{k=2}^\infty \left( \sum_{m=1}^{k-1} S_{m+1}
S_{k-m+1}
\right) z^k.
\end{equation}
Thus the equation (\ref{riccc}) becomes
\begin{eqnarray}
-\sum_{k=1}^\infty (b+k)S_{k+1}z^k 
+ \left[1-\frac{2a}{b}\right]&& \sum_{k=1}^\infty S_{k+1}z^{k+1}
+\left[\frac{a^2}{b^2} -\frac{a}{b}\right]z \nonumber \\&+&
 \sum_{k=2}^\infty \left( \sum_{m=2}^{k} S_{m}
S_{k-m+2}
\right) z^{k+1} =0   \label{comp}
\end{eqnarray}
Comparing the coefficients of $z^k,\;k=1,2,\dots$ in (\ref{comp})
we get:
$$
S_2 = \frac{a(a-b)}{b^2(b+1)},
$$
$$
S_3 = \frac{a(a-b)(b-2a)}{b^3(b+1)(b+2)},
$$
\begin{equation}
S_{k+1} = \frac{1}{b(k+b)}\left[(b-2a) S_k + b \sum_{m=2}^{k-1}
S_mS_{k-m+1} \right], \label{kishchf}
\;k=3,4,\dots.
\end{equation}
This leads, in particular, to:
$$
S_4 = \frac{a(a-b)[a(a-b)(5b+6)+b^2(b+1)]}
{b^4(b+1)^2 (b+2)(b+3)},
$$
etc., agreeing  with 
the results found by Buchholz \cite{buchholz}.

\end{document}